\documentclass[]{amsart}
\usepackage{amsmath, amssymb, amsthm,hyperref}

\begin{document}

\title{Expansive maps are isometries}

\author{Orest Bucicovschi}

\address{Department of Mathematics, UCSD, La Jolla, CA}

\email{obucicov@math.ucsd.edu}

\author{David A. Meyer}

\address{Department of Mathematics, UCSD, La Jolla, CA}

\email{dmeyer@math.ucsd.edu}

\maketitle

\begin{abstract}
In this short note we prove that any map $f$ from a dense subset $Y$ of a compact metric space $X$ into $X$ that does not decrease the distance is an isometry
\end{abstract}

\section{Preliminaries}

It is a folklore result (but see \cite{freud} ) that any map $f$ from a compact metric space $X$ to itself that does not decrease the distance ( $(d(f(x), f(y))\ge d(x,y)$)  -- also called an \emph{expansive map}, is in fact an isometry. After giving a proof of this well know statement we inquired in what measure the hypothesis $X$ complete is necessary, and we saw that in fact it is true for totally bounded metric spaces. A further modification of the proof works with maps defined on a dense subset of a totally bounded metric space. We are aware that the result may well be known, with a published but possibly different proof. Feedback would be appreciated. 

\section{Some definitions }

Let $(X,d)$ a metric space. We say that $X$ is totally bounded if for every $\epsilon>0$ there exists a covering of $Y$ with finitely many subsets of diameter $\le \epsilon$.   It is well know that a metric space is totally bounded if and only if  it is a dense subset in a compact metric space ( see \cite{kelley} ). 

Let $X$ a totally bounded metric space. Let $\epsilon >0$. There exists a covering of $X$ with $m_{\epsilon}$ subsets of diameter $\le \epsilon$.  Hence there exists at most $m_{\epsilon}$ points of $X$ with pairwise distances $>\epsilon$.  Let's define $n_{\epsilon}= n_{\epsilon}(X)$ to be the largest size of a subset of $X$  such that the distance  between any two points is $>\epsilon$.  It is easy to see that if $Y \subset X$ then $n_{\epsilon}(Y) \le n_{\epsilon}(X)$ with equality if $Y$ is dense in $X$. Let's define 
\begin{eqnarray}
\mathcal{N}_{\epsilon}(X) = \{ (x_1, \ldots, x_n) \  |\,   x_i \in X, d(x_i, x_j) > \epsilon \ \textrm{ for all } i \ne j \}
\end{eqnarray}
Let's call the elements of $\mathcal{N}_{\epsilon}(X)$ $\epsilon$-nets of $X$. It is clear that for every $\epsilon$-net ${\bf x} = (x_1, \ldots, x_{n_{\epsilon}})$ of $X$ and $y $ in $X$ there exists $i$ so that $d(y, x_i) \le \epsilon$. 

For an $\epsilon$-net ${\bf x} = (x_1, \ldots, x_{n_{\epsilon}})$ define its gauge $G({\bf x})$ to be 
\begin{eqnarray}
G({\bf x})= \prod_{i < j} d(x_i, x_j)
\end{eqnarray}
the gauge $G({\bf x})$ being a measure of the spread of ${\bf x}$. It it clear that the function $G \colon \mathcal{N}_{\epsilon}(X) \to (0, \infty)$ is bounded above  by $\textrm{diam}(X)^{n_{\epsilon}}$.  We define $g_{\epsilon}(X)$ to be the supremum of $G$ on $\mathcal{N}_{\epsilon}(X)$. Again, $Y \subset X$ implies $g_{\epsilon}(Y) \le g_{\epsilon}(X)$ with equality if $Y$ is dense in $X$. 

\section{Statement and proof of the main result}

Let $X$ a totally bounded metric space, $Y$ a dense subset of $X$ and  $f\colon Y \to X$ such that $d(f(x), f(y)) \ge d(x,y)$ for all $x$, $y$ in $Y$. Then $f$ is an isometry, that is, $d(f(x), f(y) ) = d(x,y)$ for all $x$, $y$ in $Y$. 

Proof: Let $\epsilon >0$. There exists ${\bf y} \in \mathcal{N}_{\epsilon}(Y)$ so that 
\begin{eqnarray*}
G({\bf y}) > \frac{1}{1+ \epsilon} \, \, g_{\epsilon}(X)
\end{eqnarray*}
We conclude that $f({\bf y})  \in \mathcal{N}_{\epsilon}(X)$ and moreover
\begin{eqnarray}
d(x_i, x_j) \le d(f(x_i), f(x_j)) < (1+\epsilon) d(x_i, x_j)
\end{eqnarray}

Let now $y$, $z$ be in $Y$. There exist $i$, $j$ so that 
\begin{eqnarray*}
d(f(y), f(x_i)) \le \epsilon\\
d(f(z), f(x_j)) \le \epsilon
\end{eqnarray*}
We conclude that 
\begin{eqnarray}
d(f(y), f(z)) \le d(f(x_i), f(x_j) + 2 \epsilon
\end{eqnarray}
Now, we also have 
\begin{eqnarray*}
d(y,x_i)) \le \epsilon\\
d(z,x_j)) \le \epsilon
\end{eqnarray*}
and so $d(x_i, x_j) \le d(y,z) + 2 \epsilon$.

From the above we conclude
\begin{eqnarray*}
d(f(y), f(z) ) \le d(f(x_i), f(x_j) )+ 2 \epsilon <\\
<(1+\epsilon) d(x_i, x_j) + 2 \epsilon \le (1+\epsilon) ( d(y,z) + 2\epsilon) + 2 \epsilon
\end{eqnarray*}
Since $\epsilon >0$ was arbitrary we conclude $d(f(y), f(z)) \le d(y,z)$.


\begin{thebibliography}{9}
\bibitem{kelley}
John L. Kelley 
\emph{General Topology}
D. Van Nostrand Company Inc.
 1955

 \bibitem{freud}
 Freudenthal, H., and Hurewicz, Witold. 
\emph{ Dehnungen, Verkürzungen, Isometrien}
  Fundamenta Mathematicae 
  26.1 (1936): 120-122.
 
\end{thebibliography}
\end{document}